\documentclass[reqno,centertags]{amsart}
\usepackage{amsmath,amsthm,amscd,amssymb}
\usepackage{latexsym}
\usepackage{upref}
\newcommand{\bbN}{{\mathbb{N}}}
\newcommand{\bbR}{{\mathbb{R}}}

\newcommand{\bbZ}{{\mathbb{Z}}}
\newcommand{\bbC}{{\mathbb{C}}}

\newcommand{\calB}{{\mathcal B}}

\newcommand{\calD}{{\mathcal D}}

\newcommand{\calH}{{\mathcal H}}

\newcommand{\calK}{{\mathcal K}}

\newcommand{\calN}{{\mathcal N}}

\newcommand{\calU}{{\mathcal U}}
\newcommand{\calV}{{\mathcal V}}
\newcommand{\calW}{{\mathcal W}}

\newcommand{\dott}{\,\cdot\,}
\newcommand{\no}{\nonumber}
\newcommand{\lb}{\label}

\newcommand{\ol}{\overline}

\newcommand{\loc}{\text{\rm{loc}}}

\newcommand{\ran}{\text{\rm{ran}}}

\newcommand{\bi}{\bibitem}

\renewcommand{\Re}{\text{\rm Re}}
\renewcommand{\Im}{\text{\rm Im}}

\allowdisplaybreaks

\newtheorem{theorem}{Theorem}

\newtheorem{lemma}[theorem]{Lemma}
\newtheorem{corollary}[theorem]{Corollary}

\theoremstyle{definition}
\newtheorem{definition}[theorem]{Definition}
\newtheorem{example}[theorem]{Example}

\theoremstyle{remark}

\begin{document}
\title[Krein's Formula]{An Addendum to Krein's Formula}
\author[Gesztesy, Makarov, Tsekanovskii]{Fritz Gesztesy,
Konstantin A.
Makarov, Eduard Tsekanovskii}
\address{Department of Mathematics, University of Missouri,
Columbia, MO
65211, USA}
\email{fritz@math.missouri.edu\newline
\indent{\it URL:}
http://www.math.missouri.edu/people/faculty/fgesztesy.html}
\address{Department of Mathematics, University of Missouri,
Columbia,
MO 65211, USA}
\email{makarov@azure.math.missouri.edu}
\address{Department of Mathematics, University of Missouri,
Columbia, MO
65211, USA}
\email{tsekanov@leibniz.cs.missouri.edu}
\thanks{Research supported by the US National Science Foundation
under
Grant No.~DMS-9623121.}
\date{\today}

\begin{abstract}
We provide additional results in connection with Krein's
formula, which
 describes the resolvent difference of two self-adjoint
extensions $A_1$
and $A_2$
of a densely defined closed symmetric linear operator
$\dot A$ with
deficiency
indices $(n,n),$
$n\in \bbN \cup \{ \infty \}$.~In particular, we explicitly
derive the
linear fractional transformation relating the
operator-valued Weyl-Titchmarsh
 $M$-functions $M_1(z)$ and  $M_2(z)$ corresponding to $A_1$
and $A_2$.
 \end{abstract}

\maketitle
The purpose of this note is to derive some elementary but useful
consequences of Krein's formula, which appear to have escaped
notice in
the
literature thus far.

We start with the basic setup following a short note of Saakjan
\cite{Sa65}. This paper is virtually unknown in the western
hemisphere,
and to the best of our knowledge, no English translation of it
seems to
exist. Since the paper contains few details, we provide some
proofs of
the basic facts used in \cite{Sa65}.

Let $\calH$ be a complex separable Hilbert space with scalar
product
$(\dott, \dott)$ (linear in the second factor),
denote the identity operator in $\calH$ by $I$, abbreviate the
restriction $I\big|_\calN$ of $I$ to a closed subspace
$\calN$ of
$\calH$ by $I_\calN$, and let $\calB (\calH)$
be the Banach space of
 bounded linear operators on $\calH$.
Let
$\dot A:\calD(\dot A)\to\calH, \,\, \ol{\calD(\dot A)}=\calH$
be a
densely defined closed symmetric linear operator in $\calH$ with
equal deficiency indices
 $ {\rm def}\,\, (\dot A)=(n,n),$
$n\in \bbN\cup \{ \infty \}$. We denote by $\calN_\pm$ the
deficiency
subspaces of $\dot A$,
that is,
\begin{equation}\label{e1}
\calN_\pm=\ker ({\dot A}^*\mp i).
\end{equation}
For any self-adjoint extension $A$ of $\dot A$ we introduce
its unitary
Cayley transform $C_A$ by
\begin{equation}\label{e2}
C_A=(A+i)(A-i)^{-1}.
\end{equation}

In addition, we call two self-adjoint self-adjoint extensions
$A_1$ and
$A_2$ of $\dot A$ relatively prime if
$\calD (A_1)\cap \calD (A_2)=\calD (\dot A)$. (In this case we
shall also
write that $A_1$ and $A_2$ are relatively prime w.r.t.~$\dot A$).
The point spectrum (i.e., the set of eigenvalues) and the
resolvent set
of a linear operator $T$ in $\calH$  are abbreviated by
$\sigma_p(T)$ and $\rho (T)$, respectively,
 and the direct sum of two linear subspaces $\calV$ and
$\calW$ of
$\calH$ is denoted by $\calV \dot + \calW$ in the following.

\begin{lemma} \lb{l1}
Let $A$, $A_1$, and $A_2$ be self-adjoint extensions of
$\dot A$. Then

\noindent (i).
The Cayley transform of $A$ maps $\calN_-$ onto $\calN_+$
\begin{equation}\label{e3}
C_A\calN_-=\calN_+.
\end{equation}

\noindent (ii). $\calD (A)=\calD (\dot A)
\dot + (I-C_A^{-1})\calN_+$.

\noindent (iii). For $f_+\in \calN_+$, $(A-i)^{-1}C_A^{-1}f_+
=(i/2)(C_A^{-1}-I)f_+$.

\noindent (iv). $\calN_+$ is an invariant subspace for
$C_{A_1}C_{A_2}^{-1}$ and $C_{A_2}C_{A_1}^{-1}$. In addition,
$A_1$ and $A_2$ are relatively prime if and only if
\begin{equation} \label{e4}
1\notin \sigma_p(C_{A_1}C_{A_2}^{-1}\big|_{\calN_+}).
\end{equation}

\noindent (v). Suppose $A_1$ and $A_2$ are relatively prime
w.r.t.~$\dot A$. Then
\begin{equation}\label{e5}
\ol{\ran ((A_2-i)^{-1}-(A_1-i)^{-1})}=\calN_+,
\end{equation}
\begin{equation}\label{e6}
\ker (((A_2-i)^{-1}-(A_1-i)^{-1})\big|_{\calN_-})=\{0 \}.
\end{equation}
\end{lemma}
\begin{proof}
Since the facts are standard we only sketch the main steps.\\
\noindent (i). Pick $g\in \calD(\dot A)$, $f=(\dot A-i)g$,
then
$C_Af=(\dot A+i)g \in
\ran(\dot A+i)$ yields $C_A\ran (\dot A-i)\subseteq
\ran(\dot A+i)$.
Similarly one infers
$C_A^{-1}\ran(\dot A+i)\subseteq \ran (\dot A-i)$ and hence
$C_A\ran (\dot A -i)=\ran (\dot A+i)$. Since
$C_A\in\calB(\calH)$
(in fact, $C_A$ is unitary),
$C_A\ol{\ran (\dot A-i)}=\ol{\ran (\dot A+i)}$.
$\calH=\ker (\dot A^*- i)\oplus \ol{\ran (\dot A+ i)}$ then
yields
$C_A\calN_-=\calN_+$.

\noindent (ii). By von Neumann's formula \cite{Ne29},
\begin{equation}\label{e6a}
\calD (A)= \calD (\dot A)\dot + \calN_+ \dot  +
\calU_A\calN_+
\end{equation}
for some linear isometric isomorphism $\calU_A:\calN_+
\to \calN_-$.
Since
$I-C_A^{-1}=2i(A+i)^{-1}$, $(I-C_A^{-1})\calN_+
=2i(A+i)^{-1}\calN_+
\subseteq \calD(A)$, one concludes
\begin{equation}\label{e6b}
\calU_A=-C_A^{-1}\big|_{\calN_+}.
\end{equation}

\noindent (iii). $(i/2)(I-C_A^{-1})f_+\in \calD(A)$ and
 \begin{align}
(i/2)(A-i)(I-C_A^{-1})f_+&=(i/2)(A-i)2i(A+i)^{-1}f_+
=-(A-i)(A+i)^{-1}f_+ \no \\
&=-C_A^{-1}f_+. \label{e7}
\end{align}

\noindent (iv). By (i), $\calN_+$ is an invariant subspace
for the
unitary operators $C_{A_1}C_{A_2}^{-1}$
and $C_{A_2}C_{A_1}^{-1}$. To complete the proof of
(\ref{e4})
 it suffices to note that
\begin{align}
&1\in \sigma_p(C_{A_1}C_{A_2}^{-1}\vert{\calN_+})
\Longleftrightarrow
\,\, \exists \,\, 0\ne f_+\in \calN_+ : \,\, C_{A_1}^{-1}f_+
=C_{A_2}^{-1}f_+ \no \\
&\Longleftrightarrow
(I-C_{A_1}^{-1})f_+=(I-C_{A_2}^{-1})f_+
\in \calD(A_1)\cap \calD(A_2) \lb{e7a}
\end{align}
(by the proof of (ii)) and $(I-C_{A_1}^{-1})f_+ \notin
\calD(\dot A) $
(by (ii)) yields a contradiction to $A_1$ and $A_2$ being
relatively
prime w.r.t.~$\dot A$.

\noindent  (v). Let $g\in \calD (\dot A)$, $f=(\dot A+i)g$, then
\begin{align}
&(f, ((A_2-i)^{-1}-(A_1-i)^{-1})h)
=(((A_2+i)^{-1}-(A_1+i)^{-1})(\dot A+i)g,h)=0, \no \\
&\hspace*{10.4cm} h\in \calH \label{e8}
\end{align}
yields
\begin{equation}\label{e9}
\ran ((A_2-i)^{-1}-(A_1-i)^{-1}) \subseteq \ran (\dot A+i)^\perp=
\ker (\dot A^*-i)=\calN_+.
\end{equation}
Next, let $0\ne f_+\in \calN_+$ and $f_+ \perp \ran
(( A_2-i)^{-1}
-(A_1-i)^{-1}). $
In particular,
$f_+\perp  (( A_2-i)^{-1}-(A_1-i)^{-1})C_{A_1}^{-1}f_+.$ By (iii),
$(A_1-i)^{-1}C_{A_1}^{-1}f_+=-(i/2) (I-C_{A_1}^{-1})f_+$ and
\begin{align}
(A_2-i)^{-1}C_{A_1}^{-1}f_+&=(A_2-i)^{-1}C_{A_2}^{-1}(C_{A_2}
C_{A_1}^{-1}f_+)
=-(i/2)(I-C_{A_2}^{-1})(C_{A_2}C_{A_1}^{-1}f_+) \no \\
&=-(i/2)(C_{A_2}
C_{A_1}^{-1}-C_{A_1}^{-1})f_+, \label{e10}
\end{align}
and hence
\begin{equation}\label{e11}
((A_2-i)^{-1}-(A_1-i)^{-1}) C_{A_1}^{-1}f_+=-(i/2)(C_{A_2}
C_{A_1}^{-1}-I)f_+.
\end{equation}
Thus,
$f_+ \perp (C_{A_2}C_{A_1}^{-1}-I)f_+$, that is,
\begin{equation}\label{e12}
(f_+, C_{A_2}C_{A_1}^{-1}f_+)=
\vert\vert f_+\vert\vert^2.
\end{equation}
Since $C_{A_2}C_{A_1}^{-1}\big|_{\calN_+}$ is unitary, one
concludes
$C_{A_2}C_{A_1}^{-1}f_+=f_+=C_{A_1}C_{A_2}^{-1}f_+$
and hence
\begin{equation}\label{e13}
1\in \sigma_p(C_{A_1}C_{A_2}^{-1}\vert{\calN_+}).
\end{equation}
But (\ref{e13})  contradicts the hypothesis that $A_1$ and $A_2$
are relatively prime w.r.t.~$\dot A$.
 Consequently,
$\ol{\ran ((A_2-i)^{-1}-\ran ((A_1-i)^{-1}))}=\calN_+$, which is
(\ref{e5}).

To prove (\ref{e6}) we note that every $f_-\in \calN_-$ is of
the form
$f_-=C_{A_1}^{-1}f_+$
 for some $f_+\in \calN_+$ using (i). Suppose
 $((A_2-i)^{-1}-(A_1-i)^{-1})C_{A_1}^{-1}f_+=0$. By (\ref{e11}),
this yields
$C_{A_1}C_{A_2}^{-1}f_+=f_+$ and hence
$1\in \sigma_p (C_{A_1}C_{A_2}^{-1}\big|_{\calN_+})$.
 Since $A_1$ and $A_2$ are relatively prime w.r.t.~$\dot A$ one
concludes $f_-=C_{A_1}^{-1}f_+=0$.
\end{proof}

Next, assuming $A_\ell, \ell=1,2 $ to be self-adjoint extensions
of $\dot A$, define
\begin{align}
&P_{1,2}(z)= (A_1-z)(A_1-i)^{-1} ((A_2-z)^{-1}-(A_1-z)^{-1})
(A_1-z)(A_1+i)^{-1}, \label{e14} \\
&\hspace*{8.7cm} z\in \rho(A_1)\cap\rho(A_2). \no
\end{align}
We collect the following properties of $P_{1,2}(z)$.

\begin{lemma} \mbox{\rm \cite{Sa65} } \lb{l2}
Let $z, z'\in \rho(A_1)\cap  \rho(A_2).$

\noindent (i). $P_{1,2}: \rho(A_1)\cap  \rho(A_2)\to
\calB (\calH)$
is analytic and
\begin{equation}\label{e15}
P_{1,2}(z)^*=P_{1,2}(\bar z).
\end{equation}

\noindent (ii).
\begin{equation}\label{e16}
P_{1,2}(z)\big|_{\calN_+^\perp} =0, \,\, \,\, P_{1,2}(z)\calN_+
\subseteq \calN_+.
\end{equation}

\noindent (iii).
\begin{equation}\label{e17}
P_{1,2}(z)=P_{1,2}(z')+(z-z') P_{1,2}(z')(A_1+i)(A_1-z')^{-1}
(A_1-i)(A_1-z)^{-1}P_{1,2}(z).
\end{equation}

\noindent (iv). $\ran (P_{1,2}(z)\big|_{\calN_+})$ is
independent
of  $z\in \rho(A_1)\cap\rho(A_2).$

\noindent (v). Assume $A_1$ and $A_2$ are  relatively prime
self-adjoint extensions of $\dot A$. Then
$P_{1,2}(z)\big|_{\calN_+}:\calN_+\to\calN_+ $ is invertible
(i.e., one-to-one).

\noindent (vi). Assume $A_1$ and $A_2$ are  relatively prime
self-adjoint extensions of $\dot A$. Then
\begin{equation}\label{e18}\
\ol{\ran (P_{1,2}(i))}= \calN_+.
\end{equation}

\noindent (vii).
\begin{equation}\label{e19}
P_{1,2}(i)\big|_{\calN_+} =(i/2) (I-C_{A_2}C_{A_1}^{-1})
\big|_{\calN_+}.
\end{equation}
Next, let
\begin{equation}\label{e20}
C_{A_2}C_{A_1}^{-1}\big|_{\calN_+} = -e^{-2i\alpha_{1,2}}
\end{equation}
for some self-adjoint (possibly unbounded) operator
$\alpha_{1,2} $ in $\calN_+$.
 If $A_1$ and $A_2$ are relatively prime, then
\begin{equation}\label{e21}
\{(m+\tfrac{1}{2} )\pi \}_{m\in \bbZ}
 \cap \sigma_p(\alpha_{1,2}  )=\emptyset
\end{equation}
and
\begin{equation}\label{e22}
(P_{1,2}(i)\big|_{\calN_+})^{-1} =\tan (\alpha_{1,2})
-iI_{\calN_+}.
\end{equation}
In addition,
  $\tan (\alpha_{1,2}) \in \calB (\calN_+)$ if and only if
$\ran (P_{1,2}(i)) =\calN_+$.
\end{lemma}
\begin{proof}

(i) is clear from (\ref{e14}).

\noindent (ii).  Let $f\in \calD (\dot A)$, $g=(\dot A +i)f$.
Then
\begin{equation}\label{e23}
P_{1,2}(z)g=(A_1-z)(A_1-i)^{-1}((A_2-z)^{-1}-(A_1-z)^{-1})
(\dot A-z)f=0
\end{equation}
yields $P_{1,2}(z)\big|_{\ran (\dot A+i)}=0$ and hence
$P_{1,2}(z)\big|_{\ol{\ran
(\dot A+i)}}=P_{1,2}(z)\big|_{\calN_+^\perp}=0$ since
$P_{1,2}(z) \in
\calB(\calH)$. Moreover, by
(\ref{e14})
\begin{equation}
\ran (P_{1,2}(z))\subseteq (A_1-z)(A_1-i)^{-1}
\ker (\dot A^*-z)\subseteq \ker (\dot A^*-i)=\calN_+
\end{equation}
since
\begin{align}
&(\dot A^*-i)(A_1-z)(A_1-i)^{-1}\big|_{\ker (\dot A^* -z)} \no \\
&=(\dot A^*-i)(I-(z-i)(A_1-i)^{-1})\big|_{\ker (\dot A^*-z)}
 \no \\
&=((z-i)I-(z-i)(\dot A^*-i)(A_1-i)^{-1})\big|_{\ker
(\dot A^* -z)}=0.
\label{e24}
\end{align}
This proves (\ref{e16}).

\noindent (iii). (\ref{e17}) is a straightforward (though
tedious)
computation using (\ref{e14}).

\noindent (iv). By (\ref{e17}), $\ran (P_{1,2} (z)
\big|_{\calN_+ })
\subseteq
\ran (P_{1,2}(z')\big|_{\calN_+}).$
By symmetry in $z$ and $z'$,
\begin{equation}\label{e26}
\ran (P_{1,2}(z)\big|_{\calN_+})=\ran (P_{1,2}(z')\big|_{\calN_+})
\end{equation}
is independent of $z\in \rho(A_1)\cap \rho(A_2)$.

\noindent (v). Since
\begin{equation}\label{e27}
P_{1,2}(i)=((A_2-i)^{-1} - (A_1-i)^{-1}) C_{A_1}^{-1},
\end{equation}
$C_{A_1}^{-1}: \calN_+ \to \calN_-$ is isometric,
$\ker (( A_2-i)^{-1}
-(A_1-i)^{-1}))\big|_{\calN_-})=\{0 \}$ by
(\ref{e6}), one infers $\ker ((P_{1,2}(i)\big|_{\calN_+})=\{0\}.$
Taking
$z'=i$ in (\ref{e17}) yields
$\ker (P_{1,2} (z)\big|_{\calN_+}) =\{0\}$, that is,
$P_{1,2} (z)\big|_{\calN_+}$ is invertible.

\noindent (vi). Since
$\ol{\ran ((A_2-i)^{-1}-(A_1-i)^{-1})}=\calN_+ $ by (\ref{e5}) and
$C_{A_1}^{-1}:\calH\to \calH$ is unitary, \eqref{e26}
implies (\ref{e18}).

\noindent (vii).~(\ref{e19}) follows from (\ref{e11})
and (\ref{e27}).~(\ref{e21}) is a consequence of (\ref{e4}), and
(\ref{e22}) follows
from the
 elementary trigonometric identity
$(i/2)(1+e^{-2ix})=(\tan (x)-i)^{-1}$.
By (\ref{e22})  and (\ref{e18}) $\ran (P_{1,2}(i)) =\calD
(\tan(\alpha_{1,2}))$ is dense in $\calN_+$ and hence
$\tan (\alpha_{1,2})\in \calB(\calN_+)$ if and only if
$\ran (P_{1,2}(i))$ $=\calN_+$.
\end{proof}

Next, we turn to the definition of Weyl-Titchmarsh operators
associated with self-adjoint extensions of $\dot A$.

\begin{definition} \mbox{\rm  } \lb{d3}

Let $A$ be a self-adjoint extension of $\dot A$,
$\calN\subseteq \calN_+$ a closed linear subspace of $\calN_+
=\ker (\dot A^*-i)$, and
$z\in \rho(A)$. Then the Weyl-Titchmarsh operator
$M_{A,\calN}(z) \in \calB
(\calN)$ associated with the pair $(A, \calN)$ is defined by
\end{definition}
\begin{equation}\label{e28}
M_{A,\calN}(z)=P_{\calN}(zA+I)(A-z)^{-1}P_{\calN}\big|_{\calN}=
zI_{\calN} +(1+z^2)P_{\calN} (A-z)^{-1} P_{\calN}\big|_{\calN},
\end{equation}
with $P_{\calN}$ the orthogonal projection in $\calH$ onto
$\calN$.

Weyl-Titchmarsh $m$-functions of the type \eqref{e28} have
attracted a
lot of interest since their introduction by Weyl \cite{We10}
in the
context of second-order ordinary differential operators and
their
function-theoretic study initiated by Titchmarsh \cite{Ti62}.
Subsequently, Krein introduced the concept of $Q$-functions, the
appropriate generalization of the scalar Weyl-Titchmarsh
$m$-function,
and he and his school launched a systematic investigation of
$Q$. The
literature on $Q$-functions is too extensive to be discussed
exhaustively
in this note. We refer, for instance, to \cite{Kr44}, \cite{Kr46},
\cite{KO77}, \cite{KO78}, \cite{KS66}, \cite{LT77}, \cite{TS77}
and the literature therein. Saakjan \cite{Sa65} considers a
$Q$-function
of the type \eqref{e28} in the general case where
${\rm def}(\dot A)=(n,n)$, $n\in\bbN\cup\{\infty\}$. The special
case
${\rm def}(\dot A)=(1,1)$ was also discussed by Donoghue
\cite{Do65}, who
apparently was unaware of Krein's work in this context. For a
recent
treatment of operator-valued $m$-functions we also refer to
Derkach and
Malamud \cite{DM91}, \cite{DM95} and the extensive bibliography
therein.

\begin{lemma} \mbox{\rm \cite{Sa65} } \lb{l4}
Let $A_{\ell}$, $\ell=1,2$ be relatively prime self-adjoint
extensions
of $\dot A$. Then
\begin{align}
(P_{1,2}(z)\big|_{\calN_+})^{-1}&=
(P_{1,2}(i)\big|_{\calN_+})^{-1}
-(z-i)P_{\calN_+}(A_1 +i)(A_1-z)^{-1}P_{\calN_+} \label{e29} \\
&=\tan (\alpha_{1,2})-M_{A_1, \calN_+}(z), \,\, z\in \rho(A_1).
\label{e30}
\end{align}
\end{lemma}

\begin{proof}
By Lemma~\ref{l2}\,(v), $P_{1,2}(z)\big|_{\calN_+}$ is invertible
for $z\in \rho(A_1)\cap \rho(A_2)$ and hence  (\ref{e29})
follows
from
(\ref{e27}) (and extends by continuity to all $z\in
\rho (A_1)).$
Equation (\ref{e30}) is then clear from (\ref{e22}) and
(\ref{e28}).
\end{proof}

Given Lemmas 1,2, and 4 we can summarize Saakjans's results on
Krein's formula as follows.

\begin{theorem} \mbox{\rm \cite{Sa65} } \lb{t5}
Let $A_1$ and $A_2$ be self-adjoint extensions of $\dot A$ and
$z\in \rho(A_1)\cap\rho(A_2)$. Then
\begin{align}
(A_2-z)^{-1}&=(A_1-z)^{-1}+(A_1-i)(A_1-z)^{-1}P_{1,2}(z)
(A_1+i)^{-1}(A_1-z)^{-1} \label{e31} \\
&=(A_1-z)^{-1}+(A_1-i)(A_1-z)^{-1} P_{\calN_{1,2,+}} \no \\
&\quad \, \times (\tan (\alpha_{\calN_{1,2,+}})- M_{A_1,
\calN_{1,2,+}}(z))^{-1}
P_{\calN_{1,2,+}}(A_1+i)(A_1-z)^{-1}, \label{e32}
\end{align}
where
\begin{equation} \label{e33}
\calN_{1,2,+}=\ker ((A_1\big|_{\calD (A_1)
\cap\calD (A_2)} )^* -i)
\end{equation}
and
\begin{equation} \label{e34}
e^{-2i\alpha_{\calN_{1,2,+}}}=-C_{A_2}C_{A_1}^{-1}
\big|_{\calN_{1,2,+}}.
\end{equation}
\end{theorem}
\begin{proof}
If $A_1$ and $A_2$ are relatively prime w.r.t.~$\dot A$,
Lemmas~\ref{l1}, \ref{l2}, and \ref{l4} prove
(\ref{e31})--(\ref{e34}). If $A_1$ and
$A_2$ are arbitrary self-adjoint extensions of $\dot A$ one
replaces
$\dot A$ by
the largest common symmetric part of $A_1$
and $A_2$ given by
$A_1\big|_{\calD (A_1)\cap \calD (A_2)}$.
\end{proof}

Apparently, Krein's formula (\ref{e31}), (\ref{e32}) was first
derived independently by Krein \cite{Kr44} and Naimark
\cite{Na43} in
the special case
${\rm def}(\dot A)=(1,1)$. The case ${\rm def}(\dot A)=(n,n)$,
$n\in\bbN$ is due to Krein \cite{Kr46}. A proof for this case
can also be found in the classic monograph by Akhiezer and
Glazman \cite{AG93}, Sect. 84. Saakjan \cite{Sa65} extended
Krein's formula to the general case ${\rm def}(\dot A)=(n,n)$,
$n\in\bbN\cup\{\infty\}$. In another form, the generalized
resolvent
formula for symmetric  operators (including the case of
non-densely
defined operators) has been obtained by Straus \cite{St50},
\cite{St54}.
For a variety of further results and extensions
of Krein's formula we refer, for instance, to \cite{HLS95},
\cite{HS96},
\cite{KS66}, \cite{LT77}, \cite{Ma92}, \cite{St70}, and the
literature
therein.

Saakjan  \cite{Sa65} makes no explicit attempt to relate Krein's
formula and von Neumann's parametrization
\cite{Ne29} of self-adjoint extensions of $\dot A$ (or
$A_1\big|_{\calD (A_1)\cap \calD (A_2)}$). This connection,
however,
easily follows from the preceding formalism:

\begin{corollary} \lb{c6}
\begin{equation}\label{e35}
P_{1,2}(i)\big|_{\calN_{1,2,+}}= (i/2)(I-\calU_{A_2}^{-1}
\calU_{A_1})\big|_{\calN_{1,2,+}},
\end{equation}
where
\begin{equation}\label{e36}
\calU_{A_\ell}=-C_{A_\ell}^{-1}\big|_{\calN_+}, \,\, \ell=1,2
\end{equation}
denotes the linear isometric isomorphism from
$\calN_+$ onto $\calN_-$ parametrizing the self-adjoint
extensions $A_\ell$ of $\dot A$.
\end{corollary}
\begin{proof}
Combine (\ref{e6a}), (\ref{e6b}), and (\ref{e19}).
\end{proof}

Krein's formula has been used in a large variety of
problems in mathematical physics as can be inferred from the
extensive
number of references provided, for instance, in \cite{AGHKH88}.
(A
complete bibliography on Krein's formula is beyond the scope of
this
short note.)

Next, we observe that $M_{A,\calN}(z)$ and hence
 $-(P_{1,2}(z)\big|_{\calN_+})^{-1}$ and $P_{1,2}(z)
\big|_{\calN_+}$
(cf. (\ref{e30})) are operator-valued Herglotz functions. More
precisely,
denoting $\Re(T)=(T+T^*)/2$,
$\Im (T)=(T-T^*)/2i$ for linear operators $T$ with $\calD (T)
=\calD (T^*)$ in some Hilbert space $\calK$, one can prove the
following
result.

\begin{lemma} \lb{l7}
Let $A$ be a self-adjoint extension of $\dot A$, $\calN$ a closed
subspace
of $\calN_+$. Then the Weyl-Titchmarsh operator
$M_{A,\calN}(z)$ is analytic for $z\in \bbC\backslash \bbR$ and
\begin{equation}\label{e37}
\Im(z) \Im(M_{A,\calN} (z)) \geq (\max (1,|z|^2)+|\Re(z)|)^{-1},
\quad z\in \bbC\backslash \bbR.
\end{equation}
In particular, $M_{A,\calN} (z)$ is a $\calB(\calN)$-valued
Herglotz
function.
\end{lemma}
\begin{proof}
 Using (\ref{e28}), an explicit computation yields
\begin{align}
\Im (z) \Im (M_{A,\calN} (z))&=P_{\calN} (I+A^2)^{1/2}
((A-\Re (z))^2 +(\Im(z))^2)^{-1} \no \\
& \quad \, \times (I+A^2)^{1/2}P_{\calN}\big|_{\calN}. \label{e38}
\end{align}
Next we note that for $z\in\bbC\backslash\bbR$,
\begin{equation}
\frac{1+\lambda^2}{(\lambda-\Re(z))^2+(\Im(z))^2}\geq
\frac{1}{\max (1,|z|^2)+|\Re(z)|}, \quad \lambda\in\bbR. \lb{e38a}
\end{equation}
Since by the Rayleigh-Ritz technique, projection onto a subspace
contained
in the domain of a self-adjoint operator bounded from below can
only
rise the lower bound of the spectrum (cf.
\cite{RS78}, Sect.~XIII.1), \eqref{e38} and \eqref{e38a} prove
\eqref{e37}.
\end{proof}

\vspace{3mm}

In the remainder of this note we shall explicitly derive the
linear
fractional transformation relating the  Weyl-Titchmarsh operators
$M_{A_\ell, \calN_{1,2,+}}$ associated with two self-adjoint
extensions $A_\ell, \, \ell=1,2,$ of $\dot A$. For simplicity we
first
consider the case where $A_1 $ and $A_2 $ are relatively prime
w.r.t.~$\dot A$.

\begin{lemma} \lb{l8}
Suppose $A_1$  and $A_2$ are relatively prime self-adjoint
extensions
of $\dot A$ and $z\in \rho(A_1) \cap \rho(a_2)$.
 Then
\begin{align}
M_{A_2, \calN_+}(z)&=(P_{1,2}(i)\big|_{\calN_+} +(I_{\calN_+}+i
P_{1,2}(i)\big|_{\calN_+})M_{A_1, \calN_+}(z)) \no \\
& \quad \times((I_{\calN_+}+iP_{1,2}(i)\big|_{\calN_+})-
P_{1,2}(i)\big|_{\calN_+}M_{A_1, \calN_+}(z))^{-1} \label{e39} \\
&=e^{-i\alpha_{1,2}}(\cos (\alpha_{1,2})
+\sin (\alpha_{1,2}) M_{A_1, \calN_+}(z)) \no \\
& \quad \times (\sin (\alpha_{1,2})- \cos (\alpha_{1,2})M_{A_1,
\calN_+}(z))^{-1}e^{i\alpha_{1,2}}, \label{e40}
\end{align}
where
\begin{equation}\lb{e41}
e^{-2i\alpha_{1,2}}=
-C_{A_2}C_{A_1}^{-1}\big|_{\calN_+}.
\end{equation}
\end{lemma}
\begin{proof}
Using (\ref{e30}) and (\ref{e32}) one computes
\begin{align}
&M_{A_2, \calN_+}(z)=(zI+(1+z^2)P_{\calN_+}(A_2-z)^{-1}P_{\calN_+})
\big|_{\calN_+} \no \\
&=M_{A_1, \calN_+}(z)+
(1+z^2)P_{\calN_+}
(A_1-i)(A_1-z)^{-1}
P_{\calN_+}(\tan
(\alpha_{1,2})-M_{A_1, \calN_+}(z))^{-1} \no \\
&\hspace*{7.4cm} \times P_{\calN_+}(A_1+i)(A_1-z)^{-1}P_{\calN_+}
 \no \\
&=M_{A_1, \calN_+}(z)+(iI_{\calN_+}+M_{A_1, \calN_+}(z))
(\tan(\alpha_{1,2})-M_{A_1, \calN_+}(z))^{-1} \no \\
&\hspace*{5.9cm} \times (-iI_{\calN_+}+M_{A_1, \calN_+}(z))
 \no \\
&=M_{A_1, \calN_+}(z)+(iI_{\calN_+}+M_{A_1, \calN_+}(z))
(\tan(\alpha_{1,2})-M_{A_1, \calN_+}(z))^{-1} \no \\
&\hspace*{2.2cm} \times (-iI_{\calN_+}+M_{A_1, \calN_+}(z) -
\tan(\alpha_{1,2})+\tan(\alpha_{1,2})) \no \\
&=-iI_{\calN_+}+(iI_{\calN_+}+M_{A_1, \calN_+}(z))
(\tan(\alpha_{1,2})
-M_{A_1, \calN_+}(z))^{-1}(-iI_{\calN_+}+\tan(\alpha_{1,2}))
 \no \\
&=-i(-iI_{\calN_+}+\tan(\alpha_{1,2}))^{-1}(\tan(\alpha_{1,2})
-M_{A_1, \calN_+}(z))(\tan(\alpha_{1,2})-M_{A_1, \calN_+}(z))^{-1}
 \no \\
&\hspace*{8.5cm} \times (-iI_{\calN_+}+\tan(\alpha_{1,2}) )
 \no \\
&\quad +(iI_{\calN_+}+M_{A_1, \calN_+}(z))(\tan(\alpha_{1,2})-
M_{A_1, \calN_+}(z))^{-1}(-iI_{\calN_+}+\tan(\alpha_{1,2}) )
 \no \\
&=(-iI_{\calN_+}+\tan(\alpha_{1,2}))^{-1}[-i\tan (\alpha_{1,2})+
iM_{A_1, \calN_+}(z) \no \\
&\quad +(-iI_{\calN_+}+\tan(\alpha_{1,2}))(iI_{\calN_+}+M_{A_1,
\calN_+}(z))]
(\tan (\alpha_{1,2}) -M_{A_1, \calN_+}(z))^{-1} \no \\
&\hspace{7.3cm} \times ((-iI_{\calN_+}+\tan(\alpha_{1,2})) \no \\
&=(-iI_{\calN_+}+\tan(\alpha_{1,2}))^{-1}(I_{\calN_+}
+\tan (\alpha_{1,2})M_{A_1, \calN_+}(z)) (\tan (\alpha_{1,2})
-M_{A_1,
\calN_+}(z))^{-1}\no \\
&\hspace*{6cm} \times
((-iI_{\calN_+}+\tan(\alpha_{1,2}))=(\ref{e40}). \label{e42}
\end{align}

Equation (\ref{e39}) then immediately follows
from (\ref{e42}) since
$P_{1,2}(i)\big|_{\calN_+}=(\tan (\alpha_{1,2} )
-iI_{\calN_+})^{-1}$
by (\ref{e22}).
\end{proof}

Finally, we treat the case of general self-adjoint extensions of
$\dot A$
and state the principal result of this note.

\begin{theorem} \lb{t9}
Suppose $A_1$ and $A_2$ are self-adjoint extensions
of $\dot A$ and $z\in \rho(A_1)\cap\rho(A_2)$. Then (\ref{e39})
still
holds, that is,
\begin{align}
M_{A_2, \calN_+}(z)&=(P_{1,2}(i)\big|_{\calN_+} +(I_{\calN_+}+i
P_{1,2}(i)\big|_{\calN_+})M_{A_1, \calN_+}(z)) \no \\
&\quad \times ((I_{\calN_+}+iP_{1,2}(i)\big|_{\calN_+})-
P_{1,2}(i)\big|_{\calN_+}M_{A_1, \calN_+}(z))^{-1}, \label{e43}
\end{align}
where
\begin{equation}
P_{1,2}(i)\big|_{\calN_+} =(i/2)(I-C_{A_2}C_{A_1}^{-1})
\big|_{\calN_+},
\,\, I_{\calN_+}+P_{1,2}(i)\big|_{\calN_+}=
(1/2)(I+C_{A_2}C_{A_1}^{-1})\big|_{\calN_+}. \lb{e44}
\end{equation}
\end{theorem}
\begin{proof}
 Choose a self-adjoint extension $A_3$ of $\dot A$ such that
$(A_1, A_3)$ and
$(A_2, A_3)$ are relatively prime w.r.t.~$\dot A$. (Existence of
$A_3$ is easily
confirmed using the criterion (\ref{e4})). Then express
$M_{A_1,\calN_+}(z)$ in terms of
$M_{A_3,\calN_+}(z)$ and an associated $\alpha_{3,1}$ according
to (\ref{e40}) and \eqref{e41} and similarly, express
$M_{A_2,\calN_+}(z)$ in terms of
$M_{A_3,\calN_+}(z)$ and some  $\alpha_{3,2}$. One obtains,
\begin{align}
M_{A_1, \calN_+}(z)&=e^{-i\alpha_{3,1}}
(\cos (\alpha_{3,1}) +\sin (\alpha_{3,1})M_{A_3, \calN_+}(z))
 \no \\
&\quad \times (\sin (\alpha_{3,1}) -\cos (\alpha_{3,1})M_{A_3,
\calN_+}(z))^{-1}
e^{i\alpha_{3,1}}, \lb{e45} \\
M_{A_2, \calN_+}(z)&=e^{-i\alpha_{3,2}}(\cos (\alpha_{3,2})
+\sin (\alpha_{3,2})M_{A_3,\calN_+}(z)) \no \\
&\quad \times (\sin (\alpha_{3,2}) -\cos (\alpha_{3,2})M_{A_3,
\calN_+}(z))^{-1}e^{i\alpha_{3,2}}. \lb{e46}
\end{align}
Computing
$M_{A_3, \calN_+}(z)$ from (\ref{e45}) yields
\begin{align}
M_{A_3, \calN_+}(z)&=-e^{i\alpha_{3,1}}
(\cos (\alpha_{3,1}) -\sin (\alpha_{3,1})M_{A_1, \calN_+}(z))
 \no \\
&\quad \times (\sin (\alpha_{3,1}) +\cos (\alpha_{3,1})M_{A_1,
\calN_+}(z))^{-1}
e^{-i\alpha_{3,1}}. \label{e47}
\end{align}
Insertion of (\ref{e47}) into (\ref{e46}) yields (\ref{e43})
taking into
account \eqref{e44}.
\end{proof}

Since the boundary values
$\lim_{\varepsilon\downarrow 0}(f,M_{A_1,\calN_+}
(\lambda+i\varepsilon)g)$ for $f,g\in\calN_+$ and
a.e.~$\lambda\in\bbR$
contain spectral information
on the self-adjoint extension $A_1$ of $\dot A$, relations of
the type
\eqref{e43} entail important connections between the spectra of
$A_1$ and $A_2$. In particular, the well-known unitary
equivalence
of the absolutely continuous parts $A_{1,ac}$ and $A_{2,ac}$ of
$A_1$ and $A_2$ in the case ${\rm def}(\dot A)=(n,n)$,
$n\in\bbN$, can be inferred from \eqref{e43} as discussed in
detail in \cite{GT97}. Moreover, in concrete applications to
ordinary
differential operators with matrix-valued coefficients, the
choice of
different self-adjoint boundary conditions associated with a
given
formally symmetric differential expression $\tau$ yields
self-adjoint
realizations of $\tau$ whose corresponding $M$-functions are
related via
linear fractional transformations of the type considered in
Theorem~\ref{t9}.\\

Although it appears very unlikely that the
explicit formula \eqref{e43} has been missed in the
extensive literature on self-adjoint extensions of symmetric
operators
of equal deficiency indices, we were not able to locate a
pertinent
reference. In the special case ${\rm def}(\dot A)=(1,1)$,
equations
\eqref{e40} and \eqref{e43} are of course well-known and were
studied
in great detail by Aronszajn \cite{Ar57} and Donoghue
\cite{Do65}.\\

We conclude with a simple illustration.

\begin{example}
$\calH=L^2((0,\infty);dx)$,
\begin{align}
&\dot A=-\frac{d^2}{dx^2}, \no \\
&\calD(\dot A)=\{g\in L^2((0,\infty);dx) \, \vert \,
g,g'\in AC_{\loc}((0,\infty)), \,\, g(+0)=g'(+0)=0 \}, \no \\
&\dot A^*=-\frac{d^2}{dx^2}, \no \\
&\calD (\dot A^*)=\{g\in L^2((0,\infty);dx) \, \vert \,
g,g'\in AC_{\loc}((0,\infty)),\, g''\in   L^2((0,\infty);dx) \},
 \no \\
&A_1=A_F=-\frac{d^2}{dx^2}, \quad
\calD(A_1)=\{g\in \calD(\dot A^*) \,\, \vert\, g(+0)=0 \},
 \no \\
&A_2=-\frac{d^2}{dx^2}, \quad
\calD(A_2)=\{g\in \calD(\dot A^*) \,\, \vert \,
g'(+0)+2^{-1/2}(1-
\tan (\alpha_2))g(+0)=0  \}, \no \\
& \hspace*{9.1cm} \alpha_2\in [0,\pi)\backslash\{\pi/2\}, \no
\end{align}
where $A_F$ denotes the Friedrichs extension of $\dot A$
(corresponding to
$\alpha_2=\pi/2$). One verifies,
\begin{align}
&\ker(\dot A^*-z)=\{ce^{i\sqrt{z}x}, c \in \bbC \}, \,\, \Im \,
(\sqrt{z}) >0,
\, z\in \bbC\backslash [0, \infty), \quad {\rm def} (\dot A)
=(1,1), \no \\
&(A_2-z)^{-1}=(A_1-z)^{-1}
-(2^{-1/2}(1-\tan(\alpha_2)) +i\sqrt{z})^{-1} (\ol{e^{i\sqrt{z}
 \cdot}}
, \, \cdot \, )e^{i\sqrt{z} \cdot}, \no \\
& \hspace*{7cm} z\in \rho(A_2), \,\, \Im \,(\sqrt{z})>0, \no \\
&\calU_2^{-1}\calU_1=-e^{-2i\alpha_2}, \no \\
&P_{1,2}(z)=-(1-\tan(\alpha_2) +i\sqrt{2z})^{-1}, \,\, z\in
 \rho(A_2),
\quad
P_{1,2}(i)^{-1}=\tan(\alpha_2)-i, \no \\
&M_{A_1,\calN_+}(z)=i\sqrt{2z} +1, \quad
M_{A_2,\calN_+}(z)=\frac{\cos(\alpha_2)+\sin(\alpha_2)
(i\sqrt{2z} +1)}
{\sin(\alpha_2)-\cos(\alpha_2)(i\sqrt{2z} +1)}. \no
\end{align}
\end{example}





\end{document}